\documentstyle[12pt]{article}

\catcode`\é=\active \def é{\' e}
\catcode`\ö=\active \def ö{\¨ o}

\title{Global existence and causality for a transmission problem with a repulsive nonlinearity}
\author{F. Ali Mehmeti, V. R\'egnier \thanks{Laboratoire de Mathématiques
Appliquées et de Calcul Scientifique, Institut des Sciences et Techniques de
Valenciennes, Université de Valenciennes et du Hainaut-Cambrésis, Le Mont
Houy, 59313 VALENCIENNES Cedex 9, FRANCE, courriels :
alimehme@univ-valenciennes.fr ; Virginie.Regnier@univ-valenciennes.fr}}

\date{~}

\def\R{\mathop{{\rm I}\kern-.2em{\rm R}}\nolimits}
\def\P{\hbox{I\kern-.2em\hbox{P}}}

\def\N{\hbox{I\kern-.2em\hbox{N}}}

\newcommand{\dfrac}[2]{\displaystyle \frac{#1}{#2}}

\newcommand{\cc}[1]{\cal{#1}}

\def\qqs{\ \forall \ }

\catcode`\à=\active \def à{\` a}
\catcode`\â=\active \def â{\^ a}
\catcode`\é=\active \def é{\' e}
\catcode`\è=\active \def è{\` e}
\catcode`\ê=\active \def ê{\^ e}
\catcode`\ù=\active \def ù{\` u}
\catcode`\ô=\active \def ô{\^ o}
\catcode`\ö=\active \def ö{\¨ o}
\catcode`\ç=\active \def ç{\c c}
\catcode`\î=\active \def î{\^ \i}

\newcommand{\fin}{\,\rule{1ex}{1ex}\,}
\def\C{\hbox{l\hskip-5.5ptC\/}}
\def\R{\hbox{\it I \hskip-5ptR}}
\newtheorem{theo}{Theorem}
\newtheorem{rem}{Remark}

\newtheorem{prop}{Proposition}
\newtheorem{lem}{Lemma}
\newtheorem{cor}{Corollary} 
\newtheorem{nota}{Notation}

\begin{document}
\maketitle

\noindent \bf Abstract \rm \\ It is well-known that the solution of the classical linear wave equation with compactly supported initial condition and vanishing initial velocity is also compactly supported in a set depending on time : the support of the solution at time t is causally related to that of the initially given condition. Reed and Simon have shown that for a real-valued Klein-Gordon equation with (nonlinear) right-hand side $- \lambda u^3$, causality still holds. We show the same property for a one-dimensional Klein-Gordon problem but with transmission and with a more general repulsive nonlinear right-hand side $F$. We also prove the global existence of a solution using the repulsiveness of $F$. In the particular case $F(u) = - \lambda u^3$, the problem is a physical model for a quantum particle submitted to self-interaction and to a potential step. 
\vspace{1cm}

\noindent \bf Key words \rm transmission problem, Klein-Gordon equation, repulsive nonlinearity, global existence, energy estimates, causality.\\
\\
\noindent \bf AMS \rm 35A05, 35C15, 35E15, 35L70, 35L90, 42A38.

\section{Introduction} 
\noindent For classical linear wave equations, causality is satisfied : if the initial condition has compact support $\Sigma$ and if the initial velocity vanishes, then $$\mbox{supp}(u(t, \cdot)) \subset \cc{C} \it (\Sigma,t)= \{ x \in \R / \rm \mbox{dist}\it(x, \Sigma) \leq c t \}$$ where $c$ is the wave velocity. The support stays in the light cone (cf. \cite{Br}, Section X, comment on the wave equation). \\
This well-known result has been generalized by Reed and Simon (cf. \cite{RS}) to the case of dispersive waves in three dimensions and with the nonlinear right-hand side $- \lambda u^3$. Our aim here is to prove the same properties for a transmission problem. The waves are still dispersive, the right-hand side is more general (a nonlinear function $F$) but we restrict ourselves to dimension one since we want to use the exact expression of the solution of the linearized problem (with $F \equiv 0$) which we have studied in previous papers (cf. \cite{fam2}, \cite{fam3}, \cite{Alreg2}, \cite{Alreg3}). A two-dimensional transmission problem is studied in \cite{reg} which could be analyzed with the methods developed below. \\
\\
\noindent Let us first introduce the mathematical context. We call $\Omega_1 = \Omega_2 =(0 ; +\infty)$, $I= \R^+$ and denote by $c$ and $a_k$ ($k=1 ; 2$) three positive constants with $a_2 > a_1$. 
We consider the nonlinear problem (NLP) of finding $u_k: I \times
\overline{\Omega_{k}} \rightarrow \R$ with the time variable $t \in I$ and the
space variable $x \in \Omega_k$, $k=1 ; 2$ satisfying : 

$$ \begin{array} {lllll}
(E_{k}) : & \dfrac {\partial^2 u_{k}}{\partial t^2} (t,x) - c^2
\dfrac{\partial^2 u_{k}}{\partial x^2}  (t,x) +a_k u_k(t,x) = F(u_k(t,x)), \forall (t,x)
\in \R^{+*} \times \Omega_k\\
(T_{0}): & u_1(t,0)  =  u_2 (t,0),\qqs t\in I \\
(T_{1}): & \dfrac{\partial u_1}{\partial x} (t,0^+) + \dfrac{\partial u_2}{\partial x} (t,0^+) = 0, \forall t \in I\\
(IC_{1}): & u_{k}(0,\cdot)  =  f_k, k=1 ; 2\\
(IC_{2}): & \dfrac {\partial u_{k}} {\partial t} (0,\cdot) =
0, k=1 ; 2 
\end{array} $$  

\noindent where $F$ is defined from $\R$ to $\R$ and $f_k$ from $\Omega_k$ to $\R$ for $k \in \{ 1 ; 2 \}$.\\
\\
\noindent The functions $f_k$ represent the initial condition and $F$ contains the nonlinearity of the problem. The solution formula of the linearized problem (obtained with $F \equiv 0$) is first obtained by Weidmann (cf. \cite{Weid}) and given in \cite{fam3} and \cite{Alreg2}. The meaning of the coefficients is also explained in the latter two papers~:~the propagation velocity of the waves (phase velocity) is $c$ (in the case $a_k = 0$), and the $a_k$'s are coefficients which characterize the dispersion.
$(T_0)$ and $(T_1)$ express the absence of energy loss at the central node $x=0$ 
(cf. Theorem 3, Section 3 of \cite{Alreg2}). See Section 5 of the same paper for a physical 
justification and some more details about the model. Let us just recall that the linearized problem is a model for a quantum particle submitted to a potential step. In the nonlinear problem the particular nonlinearity $F(u) = - \lambda u^3$ means that the particle is also submitted to self-interaction as well as to a potential step. Cf. Hladik and Chrysos's book about quantum mechanics for an introduction to the physical applications of our evolution problem (\cite{Hla}).  \\
\\
\noindent The first part of our work is devoted to the existence and uniqueness of global solutions for a transmission problem with two semi-infinite axes connected at one point. The problem is first expressed in an abstract way (Section 2.1) and the local existence and smoothness of solutions is stated using the theorems stated in Reed and Simon's book (cf. \cite{RS}, section X.13) which are based on the Banach fixed point theorem (cf. Section 2.2). The estimates needed to get a Lipschitz continuity property for the composition operator $f \mapsto F \circ f$ defined in Sobolev spaces only hold for sufficiently regular $F$'s (cf. Section 3, Proposition 2). Runst and Sickel's book about Sobolev spaces and Nemystkij operators has been a source of inspiration for the treament of the estimates (cf. \cite{RSi}). The investigation of the Nemytskij operator in certain functional spaces has its own interest in functional analysis, independently of the special application here (for a study in Lebesgue, Orlicz, H\"older or Sobolev spaces, see Appell, \cite{App}). Apart from being sufficiently regular, the nonlinearity $F$ is chosen to be repulsive to ensure the boundedness of the local solution and by the way its global existence. The reason is that the energy of the solution is the sum of two terms : one is the part one has even without the nonlinearity (if $F \equiv 0$) and the other one contains the nonlinear part of the energy. The first one is positive (it is directly connected to the norm of the solution in an adapted metric space) so the second one has to be positive as well (cf. 2. of the proof of Proposition 1). Otherwise, the solution could explode (its norm tends to $+ \infty$ in a finite time) since only the sum of both terms is constant. This restriction leads to what we call a "repulsive nonlinearity" which means the dispersion of the energy will increase with time. The spatial support of the solution will become larger and larger which keeps it from explosion. 
Segal was one of the first to get interested in the conditions on the nonlinearity $F$ which ensure a bounded (and so, global) solution. See \cite{S}. Haraux proved that the energy is constant with nonlinearities of the type $F(u)=g(u^2)u$ with some assumptions on $g$ which are equivalent to our conditions on $F$ in the second part of Proposition 1 (cf. \cite{Ha}). The original proof of the existence of global solutions of (NLP) with $a_1=a_2>0$ and $F(u)=- \lambda u |u|^2$ is due to J\"orgens (cf \cite{Jr}). \\
\\
\noindent In a second part, the spectral solution of the linearized problem stated in a previous work (\cite{Alreg3}) is recalled (Theorem 5 of Section 3.2) and causality is proved for the linear problem using holomorphy. An adaptation of the Paley-Wiener theorems is established which is at the core of the proof (cf. Section 3.2). Then, as in Reed and Simon's book (cf. \cite{RS}), the Banach fixed point theorem is used to get the generalization to the nonlinear problem and due to the assumption $F(0)=0$ the support is respected i.e. the support of $F \circ u$ is a subset of that of $u$. \\
Note that both parts are quite independent : global existence is not essential for causality.   \\
\\
\noindent Our interest in this particular transmission problem dates back to a few years ago. In \cite{Alreg3}, we studied the spectral solution and interpreted the phase gap between the original and reflected term as a delay in the reflection. This phenomenon is correlated to the complementary phenomenon of advanced transmission through a finite barrier for low frequency wave packets, theoretically described by Deutch and Low in 1993 (cf. \cite{D-L}). Enders and Nimtz had realized some experiments on superluminal barrier traversal before (\cite{Nim}, in 1992). Emig
has studied the same situation in three dimensions (with waveguides). His conclusion is in accordance with that of Deutch and Low : there is an apparent violation of causality. Using the integral curves of the energy flow in the space-time he has shown that the maximum of the incident wave packet and the maximum of the transmitted pulse are not related by this flow~:~there is a redistribution of the energy contained in the forward tail (\cite{Em}). \\  
Ali Mehmeti and Daikh refined the results of Deutch and Low using the same technique of restriction to frequency bands as in our paper of 2003 (\cite{Alreg3}). Here again causality is not violated by the apparent superluminal effects. At the same time Régnier extended this technique to a two-dimensional problem which had been first studied from a spectral point of view by Croc and Dermenjian (cf. \cite{CD}). \\
Transmission problems have been a subject of research for Ali Mehmeti and Nicaise since the eighties (cf. \cite{fam-3}, \cite{fam-2}, \cite{fam-1}, \cite{fam1}, \cite{nith} and \cite{ni2}). In 1994 (\cite{fam2}) Ali Mehmeti gave $L^{\infty}$-time decay estimates for the solution of our linearized problem. He found a time decay of at least $Const. t^{-1/4}$. The reduced decay rate in comparison with the full linear case ($Const. t^{-1/2}$, cf. Marshall-Strauss-Wainger \cite{MSW}) seemed to be caused by the tunnel effect.  But Mihalin\u{c}i\'{c} established in 1998 that the decay is $Const. t^{-1/2}$, being optimal for the time-space points with $\sqrt{a_i} t > x$ and $Const. t^{-m}$ with any integer $m$ holds for the other points (\cite{Mi}). \\
In \cite{Alreg1}, the authors studied the splitting of the energy flow in a star-shaped network for dispersive waves with a \it same \rm dispersion relation on each branch ($a_k = a, \qqs k \in \{1 \cdots n \}$). The splitting
of the energy flow is independent of the frequency in this case. Therefore the phenomenon of delay cannot appear.\\
In the future one could study the energy flow for diffraction problems, since Laplace-type solution formulae exist (cf. Ali Mehmeti \cite{fam3} and Rottbrand \cite{Rott}), as well as for scattering problems. \\ 
Concerning nonlinear evolution problems, Haraux's book gives the techniques for proving global existence and behaviour of solutions (\cite{Ha}). Some results about global existence of solutions have been established without assuming the nonlinearity to be repulsive (cf. \cite{RS}, \cite{Kl}, \cite{St} and \cite{Sh}) : dispersion is proved to make up for attractivity using an $L^{\infty}$-time decay estimate for the solution (cf. \cite{MSW}).  \\
As for Kuksin (\cite{K}), he considers a different notion of flow of energy, which is the phenomenon of the energy transition to higher frequencies for a nonlinear wave equation under n-dimensional periodic boundary conditions in a Sobolev phase space. Note that the nonlinearity keeps the energy from staying in a frequency band.

\section{Existence and uniqueness of solutions}

\subsection{Abstract setting}
\noindent Let us consider the data as in the introduction and define the linear operator $(A,D(A))$ as follows :\\
$D(A)=\{v\in\prod_{k=1}^2H^2(\Omega_k)/ v
\hspace{2mm} \mbox{satisfies} \hspace{2mm} (T_0),(T_1)\}$
\\
$A:v\in D(A)\mapsto Av:= \left(-c^2\dfrac{d^2v_k}{dx^2} + a_k v_k 
\right)_{k\in\{1 ; 2 \}}$ in $H:=~\prod_{k=1}^2 L^2(\Omega_k)$\\ and given 
$u(t)=(u_1(t,\cdot),u_2(t,\cdot))$ and $F(u(t))=(F(u_1(t,\cdot)), F(u_2(t,\cdot)))$, the initial boundary value problem (NLP) (see the introduction) can be rewritten as the following problem denoted by (P) :

$$(P) \left \{ \begin{array} {lll}
u \in C^2(\R,H), u(t)\in D(A), \qqs t\in \R^+\\
\dfrac{d^2u}{dt^2}(t)+Au(t)=F(u(t)) \quad \mbox{in} \quad H, \qqs t\in \R^+\\\\
u(0)= \Phi, \dfrac{du}{dt}(0)=0
\end{array} \right.$$

\noindent where $\Phi=(f_1, f_2) \in D(A)$ with $f_k : \Omega_k \rightarrow \R$ for $k \in \{ 1 ; 2 \}$. \\ 
\noindent Note that for the sake of simplicity, the initial velocity $u_t(0)$ is chosen to vanish but putting $u_t(0)=\Psi$ would not change the problem fundamentally. \\ 
The spectral theory of the linearized problem ($F \equiv 0$) is studied in
\cite{fam3}, in Theorem 1.5.1. Here $d_j=1, \qqs j \in \{ 1 \cdots n \}$
and for us the $c_j$'s and $a_j$'s have to be interchanged. The operator $A$ 
is proved to be self-adjoint and the existence and uniqueness of a solution of 
the linearized problem (LP) ((LP) is (P) with $F \equiv 0$) in $C^2(\R^+, \prod_{k=1}^n 
L^2(\Omega_k)) \cap C^1(\R^+, D(\sqrt{A})) \cap C^0(\R^+, D(A))$ is stated in \cite{fam3} (in the context of the analysis of the asymptotic time behaviour of the solution). Recall that
$D(\sqrt{A}):=\{v\in\prod_{k=1}^2H^1(\Omega_k)/ v \hspace{2mm} \mbox{satisfies} \hspace{2mm} (T_0)\}$.

\noindent As Reed and Simon do in \cite{RS}, one can rewrite (P) as (Q) :
$$(Q) \left \{ \begin{array} {lll}
\dfrac{d \varphi}{dt}(t)+ i \cc A\it \varphi(t)= J(\varphi(t)) \quad \mbox{in} \quad H\\
\varphi(0)= \varphi_0
\end{array} \right.$$

\noindent where $\varphi(t)$= $\left ( \begin{array} {ll}
u(t)\\
\dfrac{\partial u}{\partial t}
\end{array} \right)$ and $J(\varphi(t))$= $\left ( \begin{array} {ll}
\hspace{0.5cm} 0\\
F(u(t))
\end{array} \right)$ and $\cc{A}$= $i \left ( \begin{array} {ll}
\hspace{0.4cm} 0 \hspace{1cm} I\\
-A \hspace{1cm} 0 
\end{array} \right)$ 

\noindent The operator $\cc{A}$ is proved to be self-adjoint and closed on its domain of definition $D(A) \oplus D(\sqrt{A})$ (cf. the proposition of Section X.13 in \cite{RS}).

\subsection{Reed-Simon's Theorems for the existence of solutions}

\noindent Let us recall the three main theorems given in \cite{RS}. Reed and Simon's proofs are all based on the contraction mapping principle applied to the mapping contained in the reformulation of Problem (P) as the integral equation 
$$\varphi(t) = e^{-i \cc A \it t} \varphi_0 + \int_0^t e^{-i (t-s) \cc A \it} J(\varphi(s)) ds  $$

\noindent For any $\varphi \in D(\cc A)$, the norm of $\varphi = (u,v)$ is defined by Reed and Simon as $\| \varphi \| = (\| \sqrt{A} u\|_2^2 + \| v \|_2^2)^{1/2}$. \\
Thus, in our case, $\cc A \varphi$, $J(\varphi)$ and $\cc A \it J(\varphi)$ are defined as 
$$ \left \{ \begin{array} {lll}
\| \cc A \varphi \it \| = (\| \sqrt{A} v\|_2^2 + \| A u \|_2^2)^{1/2} \\
\| J(\varphi) \| = \| F \circ u \|_2\\
\| \cc A \it J(\varphi) \| = \| \sqrt{A} (F \circ u) \|_2
\end{array} \right.$$ 
\noindent Note that, in our case, $u(t)$ has two components so $\| u \|_2$ has to be understood as 
$$\left(  \sum_{k=1}^2 \int_{\Omega_k} |u_k(t,x)|^2 dx  \right)^{1/2}$$ 

\noindent In the following three subsections, we recall general results on the existence and smoothness of the solutions of $(Q)$, where $\cc A$ is any self-adjoint operator and $J$ any mapping from $D(\cc A)$ to $D(\cc A)$ (see in \cite{RS} for the proofs).

\subsubsection{Local existence}

\begin{theo} ~~\\
Let $\cc{A}$ \it be a self-adjoint operator on a Hilbert space $H$ \it and J a mapping from $D(\cc{A} \it)$ to $D(\cc{A} \it)$ which satisfies $J(0) = 0$ and :
  
$$\begin{array} {ll}
(H_0^L) \qquad \| J(\varphi) - J(\psi) \| \leq C(\| \varphi \|, \| \psi \|) \cdot \|\varphi - \psi \| \\
(H_1^L) \qquad \| \cc{A} \it (J(\varphi) - J(\psi)) \| \leq C(\| \varphi \|, \| \cc{A} \it \varphi \|, \| \psi \|, \| \cc{A} \psi \it \|) \cdot \| \cc{A} \it \varphi - \cc{A} \it \psi \|  
\end{array}$$

\noindent for all $(\varphi, \psi) \in D(\cc{A})\it ^2$ where each constant $C$ is a monotone increasing (everywhere finite) function of the norms indicated. Then, for each $\varphi_0 \in D(\cc{A} \it )$, there is a $T > 0$ so that Problem $(Q)$ of Section 2.1 has a unique continuously differentiable solution for $t \in [0 ; T)$. For each set of the form $\{ \varphi \
\hspace{0.05 cm} / \hspace{0.05 cm} \| \varphi \| \leq a, \| \cc{A} \it \varphi \| \leq b \}$, $T$ can be chosen uniformly for all $\varphi_0$ in the set.
\end{theo}

\begin{rem}~~\\
\noindent Note that the assumptions $(H_0)$ and $(H_1)$ of Reed and Simon (see \cite{RS}) are implied by $J(0)=0$, $(H_0^L) $ and $(H_1^L)$. Moreover $(H_0)$ and $(H_1)$ make the condition $J(0)=0$ necessary. Hence the simplifications we have made here.   
\end{rem}

\subsubsection{Local smoothness}

\begin{theo} ~~\\
(a) Let $\cc{A}$ \it be a self-adjoint operator on a Hilbert space $H$ \it and $n$ a positive integer. Let J be a mapping which takes $D(\cc{A} \it^j)$ to $D(\cc{A} \it^j)$ for all $1 \leq j \leq n$ (for all $0 \leq j \leq n$) and which satisfies $J(0)=0$ and :
  
$$\begin{array} {ll}
(H_j^L) \qquad \| \cc{A} \it^j (J(\varphi) - J(\psi)) \| \leq C(\| \varphi \|, \| \psi \|,..., \| \cc{A} \it^j \varphi \|, \| \cc{A} \it^j \psi \it \|) \cdot \| \cc{A} \it^j \varphi - \cc{A} \it^j \psi \|
\end{array}$$

\noindent for all $(\varphi, \psi) \in D(\cc{A} \it^j)^2$ where each constant $C$ is a monotone increasing (everywhere finite) function of all its arguments. Then, for each $\varphi_0 \in D(\cc{A} \it^n)$, $n \geq 1$, there is a $T_n$ so that Problem $(Q)$ of Section 2.1 has a unique solution $\varphi(t)$ for $t \in [0,T_n)$ with $\varphi(t) \in D(\cc A \it^n)$ for all $t \in [0, T_n)$. For each set of the form $\{ \varphi \ \hspace{0.05 cm} / \hspace{0.05 cm} \| \cc{A} \it^j \varphi \| \leq a_j ; j = 0,...,n \}$, $T_n$ can be chosen uniformly for all $\varphi_0$ in the set.\\
\\
\noindent (b) In addition to the hypotheses in (a), assume that for each $j < n$, $J$ has the following property : if a solution $\varphi$ is $j$ times strongly continuously differentiable with $\varphi^{(k)}(t) \in D(\cc{A} \it^{n-k})$ and $\cc{A} \it^{n-k} \varphi^{(k)}(t)$ is continuous for all $k \leq j$, then $J(\varphi(t))$ is $j$ times differentiable, $d^j J(\varphi(t))/dt^j \in D(\cc{A} \it^{n-j-1})$, and $\cc{A} \it^{n-j-1} d^j J(\varphi(t))/dt^j$ is continuous. Then the solution given in part (a) is $n$ times strongly differentiable in $t$ and $d^j J(\varphi(t))/dt^j \in D(\cc{A} \it^{n-j-1}) $.  
\end{theo}

\subsubsection{Global existence}

\begin{theo} ~~\\
Let $\cc{A}$ \it be a self-adjoint operator on a Hilbert space $H$ \it and $n$ a positive integer. Let J be a mapping which takes $D(\cc{A} \it^j)$ to $D(\cc{A} \it^j)$ for all $1 \leq j \leq n$ (for all $0 \leq j \leq n$) and which satisfies $J(0)=0$ and :
  
$$\begin{array} {lll}
(H^{\prime}_j) \qquad \| \cc{A} \it^j J(\varphi) \| \leq C(\| \varphi \|,..., \| \cc{A}\it^{j-1} \varphi \|) \cdot \| \cc{A} \it^j \varphi \|, 1 \leq j \leq n \\
(H_j^L) \qquad \| \cc{A} \it^j (J(\varphi) - J(\psi)) \| \leq C(\| \varphi \|, \| \psi \|,..., \| \cc{A} \it^j \varphi \|, \| \cc{A} \it^j \psi \it \|) \cdot \| \cc{A} \it^j \varphi - \cc{A} \it^j \psi \|\\
0 \leq j \leq n 
\end{array}$$

\noindent for all $(\varphi, \psi) \in D(\cc{A} \it^j)^2$ where each constant $C$ is a monotone increasing (everywhere finite) function of all its arguments. Let $\varphi_0 \in D(\cc{A} \it^n)$ and suppose that, on any finite interval of existence, the solution $\varphi(t)$ guaranteed by part (a) of Theorem 2 has the property that $\| \varphi(t) \|$ is bounded from above. Then there is a strongly differentiable $D(\cc A \it^n)$-valued function $\varphi(t)$ on $[0 ; + \infty)$ that satisfies 

$$\left \{ \begin{array} {lll}
\dfrac{d \varphi}{dt}(t)+ \cc A\it \varphi(t)= J(\varphi(t)) \quad \mbox{in} \quad H\\
\varphi(0)= \varphi_0
\end{array} \right.$$

\noindent Further, if J satisfies the hypotheses of part (b) of Theorem 2, then $\varphi(t)$ is $n$-times strongly differentiable and $d^j \varphi(t)/dt^j \in D(\cc A \it^{n-j})$.  
\\
\end{theo}

\subsection{Required estimates}

\noindent In this section, we will apply the last three theorems to our particular operator $\cc A$ and mapping $J$ that is we will first interpret the required estimates of Theorem 3 in terms of a Lipschitz-type continuity for the nonlinearity $F$ (Proposition 1) and establish the required regularity properties for $F$ to satisfy those estimates (Proposition 2). At last we will deduce the assumptions on $F$ necessary for the existence of a global solution to Problem $(P)$ (Corollary 1). Analogous results can be found in Haraux's book (A. II. Lecture 3 of \cite{Ha}) : in particular he proves that the energy is constant for nonlinear wave equations with a nonlinearity of the type $F(u) = g(u^2)u$ under some assumptions on $g$ in dimension $N$ ($g'$ has to lie in $L^{\infty}_{loc}(\R^+)$ in dimension one). Note that we restrict ourselves to the proof of the one-time strong differentiability of the solution ($n=1$ in Theorem 3). To get a higher regularity, nonlinear compatibility conditions have to be added to the transmission conditions $(T_0)$ and $(T_1)$ (cf. \cite{fam-1}).

\begin{prop}~~\\

\begin{enumerate}
\item For the conditions $J(0)=0$, $(H_0^L)$ and $(H_1^L)$ to be satisfied, it is 

\begin{itemize}
\item necessary that $F(0)=0$
\item sufficient that 
$$\left \{ \begin{array}{lll}
F(0)=0 \\
\| F \circ u - F \circ v \|_2 \leq C(\| u \|_2, \| u_x \|_2, \| v \|_2, \| v_x \|_2) \cdot \| \sqrt{A} u - \sqrt{A} v \|_2 \\
\| \sqrt{A} (F \circ u - F \circ v) \|_2 \leq C(\| u \|_2, \| u_x \|_2, \| v \|_2, \| v_x \|_2) \cdot \| A(u - v) \|_2
\end{array} \right.$$

\noindent for all $(u, v) \in D(A)^2$ and where $C$ is a monotone increasing function of all its arguments.

\end{itemize} 

\item Assume $F \in C^0(\R)$ and suppose that for any $(u,v) \in (\prod_{k=1}^2 L^2(\Omega_k))^2$, there exists $C(\| u\|, \| v \|)$, monotone increasing function of all its arguments such that 
$$\| F \circ u - F \circ v \|_2 \leq C(\| u \|_2, \| v \|_2) \cdot \| u - v \|_2$$ 
If $G$, defined on $\R$ by $G(w) = \int_0^w F(\xi) d \xi$, takes only non positive values and if it is such that $G \circ u$ belongs to $L^1(\Omega)$ for $u \in H^2(\Omega)$ ($\Omega$ being an interval of $\R$, possibly unbounded), then the unique local solution of Problem (P), guaranteed by the first part of the proposition in combination with Theorem 1 and denoted by $\varphi$, has the property that $\| \varphi(t) \|$ is bounded from above on any finite interval of existence. Thus, due to Theorem 3, it defines a global solution of (P) i.e. a strongly differentiable $D(\cc A)$-valued function on $[0 ; + \infty)$ that satisfies Problem (P).   

\end{enumerate}
\end{prop}

\noindent \bf Proof. \rm \\
Let us recall that the norm of $\varphi = (u,v)$ is defined as $\| \varphi \| = (\| \sqrt{A} u\|_2^2 + \| v \|_2^2)^{1/2}$ for any $\varphi \in D(\cc A)$. And those of $\cc A \varphi$, $J(\varphi)$ and $\cc A \it J(\varphi)$ are defined as 
$$ \left \{ \begin{array} {lll}
\| \cc A \varphi \it \| = (\| \sqrt{A} v\|_2^2 + \| A u \|_2^2)^{1/2} \\
\| J(\varphi) \| = \| F \circ u \|_2\\
\| \cc A \it J(\varphi) \| = \| \sqrt{A} (F \circ u) \|_2
\end{array} \right.$$

\begin{enumerate}

\item 
\begin{itemize}
\item The necessary condition is clearly required : $J(0) = 0$ and $J(\varphi(t))$= $\left ( \begin{array} {ll}
\hspace{0.5cm} 0\\
F(u(t))
\end{array} \right)$ imply $F(0)=0$.   

\item Now, $\varphi$ is still defined as $\varphi = (u,v)$ and $\psi$ as $\psi = (\tilde{u},\tilde{v})$. So, if 
$$\| F \circ u - F \circ \tilde{u} \|_2 \leq C(\| u \|_2, \| u_x \|_2, \| \tilde{u} \|_2, \| \tilde{u}_x \|_2) \cdot \| \sqrt{A} u - \sqrt{A} \tilde{u} \|_2$$
\noindent then $\| F \circ u - F \circ \tilde{u} \|_2$
$$ \leq C(\| u \|_2, \| u_x \|_2, \| \tilde{u} \|_2, \| \tilde{u}_x \|_2) \cdot (\| \sqrt{A} u - \sqrt{A} \tilde{u} \|_2^2 + \| v - \tilde{v} \|_2^2)^{1/2}$$  
\noindent Now, the norm of $\sqrt{A} u$ is defined by :
$$\| \sqrt{A} u \|_2 = (\sqrt{A} u ; \sqrt{A} u) = (Au, u) = \sum_{k=1}^2 \| u_k \|^2 + a_k \| u_{k,x} \|^2$$     
so $\| u_x \|_2 \leq \max(1/a_1 ; 1/ a_2) \| \sqrt{A} u \|_2 \leq \max(1/a_1 ; 1/ a_2) \| \varphi \|$ and since $C$ is a monotone increasing function of all its arguments, $C(\| u \|_2, \| u_x \|_2)$ is lower than $C^{\prime} ( \| \varphi \|)$, with $C^{\prime} ( \| \varphi \|)$ a monotone increasing function that will be denoted by $C$ as well. \\
Likewise $C(\| u \|_2, \| u_x \|_2, \| \tilde{u} \|_2, \| \tilde{u}_x \|_2) \leq C( \| \varphi \|, \| \psi \|)$  \\
and $(\| \sqrt{A} u - \sqrt{A} \tilde{u} \|_2^2 + \| v - \tilde{v} \|_2^2)^{1/2} = \| \varphi - \psi \|$. \\
So the first condition given in the first part of the theorem, has been proved to be sufficient for $(H_0^L)$ to hold. Idem for $H_1^L$.

\end{itemize}

\item If $u$ is the solution of Problem $(P)$, the energy is defined as
$$E(t) = \dfrac{1}{2} \sum_{k=1}^2 \int_{\Omega_k} \left( \left|(\sqrt{A} u(t,x))_k \right|^2 + |u_{k,t}(t,x)|^2 - 2 \hspace{0.2cm} (G \circ u_k)(t,x) \right) dx$$
\noindent that is : $E(t) = \dfrac{1}{2} \| \varphi(t)  \|^2 - \sum_{k=1}^2 \int_{\Omega_k} (G \circ u_k)(t,x) dx$. \\
Note that $G \circ u_k$ belongs to $L^1(\Omega_k)$ by assumption ($u$ is the solution so $u_k(t,\cdot)$ lies in $H^2(\Omega_k)$). Thus, since $G$ takes only non positive values,  
$$\dfrac{1}{2} \| \varphi(t) \|^2  \leq \dfrac{1}{2} \| \varphi(t)  \|^2 - \sum_{k=1}^2 \int_{\Omega_k} (G \circ u_k)(t,x) dx$$ i.e. $\dfrac{1}{2} \| \varphi(t) \|^2  \leq E(t)$. \\
Now let us show that the energy is independent of $t$ to get $\dfrac{1}{2} \| \varphi(t) \|^2  \leq E(0)$ i.e. a bounded norm of $\varphi(t)$ on its interval of existence. \\
The way it will be done is the same one as in Haraux (\cite{Ha},~p.~23). His nonlinearity is of the form $F(u)=u \cdot g(u^2)$ and his conditions on $g$ in dimension one are equivalent to ours. \\
So first of all, the first two terms in $E(t)$ are clearly differentiable. The reason is the same one as in Reed and Simon's book : $\varphi$ is strongly differentiable as a $D(\cc{A} \it)$-valued function. It means in particular that :
$$\left \| \sqrt{A} \left( \frac{u(t+h)-u(t)}{h} - u_t(t) \right)  \right \|_2 \rightarrow 0$$     
$$\left \| \left( \frac{u_t(t+h)-u_t(t)}{h} - u_{tt}(t) \right)  \right \|_2 \rightarrow 0$$     
\noindent as $h \rightarrow 0$. \\
There remains to deal with the third term containing the nonlinearity : \\
Since $F \in C^0(\R)$ and since $G$ is defined on $\R$ by $G(w) = \int_0^w F(\xi) d \xi$, $G$ is differentiable and $G' = F$. Now, for any $(w, l) \in \R^2$, there exists $c$ in $[\min(w;w+l) ; \max(w;w+l)]$ such that, 

$$G(w + l) - G(w) = l F(c)$$

\noindent Applied to $w = u_k(t,x)$ and $l = u_k(t+h,x)- u_k(t,x)$, it reads : 
$$G(u_k(t+h,x))-G(u_k(t,x))= F(y_k^h(t,x)) \cdot (u_k(t+h,x)-u_k(t,x))$$
\noindent with $y_k^h(t,x)  \in [\min(u_k(t,x); u_k(t+h,x)) ; \max(u_k(t,x); u_k(t+h,x))]$. \\
Note that we should write $F(y^h(t,x))_k$ with $y^h(t,x) = (y_1^h(t,x) ; y_2^h(t,x))$ since $F$ has been defined on the product space $H$. We will always keep this notation in the following. Now  

$$\begin{array}{lll}
\int_{\Omega_k} \left| \frac{G(u_k(t+h,x))-G(u_k(t,x))}{h} - F(u_k(t,x)) u_{k,t}(t,x) \right| dx \leq \\
\\
\int_{\Omega_k} \left| \frac{1}{h} \left \{ [G(u_k(t+h,x))-G(u_k(t,x))] - F(u_k(t,x)) (u_k(t+h,x)-u_k(t,x)) \right \} \right|dx \\
\\
+ \int_{\Omega_k} \left| F(u_k(t,x)) \left( \frac{u_k(t+h,x)-u_k(t,x)}{h}- u_{k,t}(t,x) \right) \right| dx 
\end{array}$$    

\noindent Then, due to H\"older's inequality
$$\begin{array}{llll} 
\sum_{k=1}^2 \int_{\Omega_k} \left| \frac{G(u_k(t+h,x))-G(u_k(t,x))}{h} - F(u_k(t,x)) u_{k,t}(t,x) \right| dx  \\
\leq \left( \sum_k \| F \circ y_k^h(t) - F \circ u_k(t) \|_2 \right) \left \| \frac{u(t+h)-u(t)}{h} \right \|_2 + \| F(u(t)) \|_2 \left \| \frac{u(t+h)-u(t)}{h} - u_t(t) \right \|_2 \\
\leq \left( \sum_k C(\| y_k^h(t) \|_2 ; \| u_k(t) \|_2) \cdot \| y_k^h(t) - u_k(t) \|_2 \right) \left \| \frac{u(t+h)-u(t)}{h} \right \|_2 \\
\hspace{0.5cm} + \| F(u(t)) \|_2 \left \| \frac{u(t+h)-u(t)}{h} - u_t(t) \right \|_2 
\end{array}$$    

\noindent Now $\frac{u(t+h)-u(t)}{h}$ tends to $u_t(t)$ in $L^2$ since $u \in C^1(\R,H)$ (it is even of class $C^2$). So there remains to prove that $\| F \circ y_k^h(t) - F \circ u_k(t) \|_2$ tends to zero with $h$ to conclude that the latter estimate tends to 0, which means that the third term is differentiable and the derivative of $E(t)$ is then  
$$E'(t) = \sum_{k=1}^2 \int_{\Omega_k} u_{k,t}(t,x) \left( u_{k,tt}(t,x)+(Au)_k(t,x) - [F(u(t,x))]_k \right) = 0$$ 
\noindent due to Equation $(E_k)$. \\
\\
First of all, it holds classically :
$$\left \{ \begin{array}{ll}
\min(l;l') = \dfrac{1}{2} (l + l' - |l - l'|)\\
\max(l;l') = \dfrac{1}{2} (l + l' + |l - l'|)\\
\end{array} \right.$$
Then 
$$\begin{array}{ll}
|y_k^h(t,x)|^2 \leq \dfrac{1}{2} (u(t,x) + u(t+h,x) + |u(t,x) - u(t+h,x)|)^2 \\
\hspace{1.8cm} \leq 4 (|u(t,x)|^2 + |u(t+h,x)|^2)
\end{array}$$

\noindent and $y_k^h(t, \cdot)$ belongs to $L^2(\Omega_k)$ since both $u(t, \cdot)$ and $u(t+h, \cdot)$ are in $L^2$. Moreover $\| y_k^h(t) \|_2^2 \leq  4 (\|u(t) \|_2^2 + \| u(t+h) \|_2^2)$.\\
Now $y_k^h(t,x)$ is superior to $m := \dfrac{1}{2} (u(t,x) + u(t+h,x) - |u(t,x) - u(t+h,x)|)$ and inferior to $M := \dfrac{1}{2} (u(t,x) + u(t+h,x) + |u(t,x) - u(t+h,x)|)$.

\noindent Then $|y_k^h(t,x) - u(t,x)|$ is inferior to the maximum of $m$ and $M$. So it is bounded by
$\left( u(t+h,x) - u(t,x) \right)$. Thus $\| y_k^h(t) - u(t) \|_2 \leq \| u(t+h) - u(t) \|_2$ and, since $u \in C^0(\R,H)$, $u(t+h)$ tends to $u(t)$ in $L^2$. \\
Then $\lim_{h \longrightarrow 0} \| y_k^h(t) - u(t) \|_2 = 0$.  \fin 

\end{enumerate}

\noindent Conditions on the nonlinearity of the problem i.e. on the function $F$ can be deduced from the latter proposition.  Runst and Sickel's book (cf. \cite{RSi}) has been a source of inspiration for the following and especially for the next two proofs. Using the Fatou property of the Sobolev space $H^2((0;+\infty))$ is their key idea (cf. \cite{RSi}, p. 15).

\begin{nota} ~~\\
The composition operator $T_F : 
\begin{array} {ll}
B \longrightarrow B'\\
f \mapsto F \circ f
\end{array}$
\noindent satisfies Property $\cc{L} \it ip\it(B,B')$ if there exists a constant $C > 0$, monotone increasing function of $\|f \|_2$, $\|g \|_2$, $\|f_x \|_2$, $\|g_x \|_2$, such that for any $(f,g)$ in $B^2 \subset (L^2)^2$
$$\|F \circ f - F \circ g \|_{B'} \leq C(\|f \|_2,\|g \|_2, \|f_x \|_2, \|g_x \|_2) \cdot \|f - g \|_B$$
\noindent where $\| \cdot \|_2$ is the norm in $L^2$.
\end{nota}

\begin{prop} ~~\\

\noindent Let $F : \R \longrightarrow \R$ be a continuous function satisfying $F(0)=0$. Denote by $T_F$ the composition operator $T_F : f \mapsto F \circ f$

\noindent Then, for $T_F$ to satisfy Property $\cc{L} \it ip\it(H^1((0;+\infty)), L^2((0;+\infty)))$ (resp. $\cc{L} \it ip\it(H^2, H^1))$, it is :   
\begin{itemize}
\item necessary that $F \in H^{2,loc}((0;+\infty))$ \\
(respectively $F \in H^{1,loc}((0;+\infty))$) 
\item sufficient that $F \in C^2(\R,\R)$ (respectively $F \in C^1(\R,\R)$) 
\end{itemize} 

\end{prop}

\noindent \bf Proof. \rm \\

\noindent Runst and Sickel's proof is adapted since the estimate we need is a weaker result as theirs : it contains a constant depending on four parameters, so we will explicit this function and check its monotony. 
\begin{enumerate}
\item Necessary condition : Choosing a sequence of functions $f_M$ in $\cc{S} \it(\R;\R)$ such that $\mapsto f_M(x)=x \quad \mbox{when} \quad |x| \leq M$ leads to $T_F(f_M)(x) = (F \circ f_M)(x) = F(x) \quad \mbox{when} \quad |x| \leq M$ i.e. $T_F(f_M) = F|_{[-M;M]}$. Thus for $T_F(f_M)$ to belong to $H^2((0;+\infty))$ (respectively $H^1((0;+\infty))$), $F|_{[-M;M]}$ has to lie in $H^2((0;+\infty))$ (resp. $H^1((0;+\infty))$) for any $M$ in $\R$ i.e. $F$ has to lie in $H^{2,loc}((0;+\infty))$ (resp. $H^{1,loc}((0;+\infty))$). 

\item Sufficient condition : the aim is to prove that, if $F \in C^2(\R)$, $\forall (f,g)$ in $(H^2((0;+\infty)))^2, \exists C(\|f \|,\|g \|, \|f_x \|, \|g_x \|)  > 0$ \\
$$\|F \circ f - F \circ g \|_{H^1}\leq C(\|f \|,\|g \|, \|f_x \|, \|g_x \|)  \|f - g \|_{H^2}$$
\noindent i.e. $T_F$ satisfies Property $\cc{L} \it ip\it(H^2, H^1)$. 
\noindent Since the arguments for $\cc{L} \it ip\it(H^1, L^2)$ are strictly analogous, this case is not treated separately.

\begin{itemize}

\item \it First step. \rm \\
Let us first prove that Property $\cc{L} \it ip \it(C_c^{\infty}(\R) \cap H^2((0;+\infty)), H^1((0;+\infty)))$ is satisfied for $F \in C^2(\R)$. \\
\noindent The norm $\| F \circ f - F \circ g \|_{H^1(\Omega)}$ is equivalent to $\| F \circ f - F \circ g \|_{L^2(\Omega)} + \| \partial_x(F \circ f - F \circ g) \|_{L^2(\Omega)}$ for any domain $\Omega$. So the aim is to find estimates for $\| F \circ f - F \circ g \|_{L^2(\Omega)}$ and $\| (F' \circ f)(f') - (F' \circ g)(g') \|_{L^2(\Omega)}$ assuming that $(f,g) \in (C_c^{\infty}(\R) \cap H^2((0;+\infty)))^2$ and $F \in C^2(\R)$.    \\

\begin{enumerate}
\item $\| F \circ f - F \circ g \|_{L^2((0;+ \infty))}^2 = \int_0^{+ \infty} |(F \circ f)(x) - (F \circ g)(x)|^2 dx$ and $|(F \circ f)(x)- (F \circ g)(x)| = |F(f(x))-F(g(x))|$. Thus, for any $x \in (0 ; + \infty)$ :
$$|((F \circ f)- (F \circ g))(x)| \leq \left[ \sup_{|y| \leq max(\| f \|_{L^{\infty}} ; \| g \|_{L^{\infty}})} |F'(y)| \right] \cdot |f(x) - g(x)|$$
\noindent So $\| F \circ f - F \circ g \|_{L^2((0;+ \infty))} \leq M(F,f,g) \cdot \| f - g \|_{L^2((0;+ \infty))}$ with 
$$M(F,f,g) = \sup_{|y| \leq max(\| f \|_{L^{\infty}} ; \| g \|_{L^{\infty}})} |F'(y)|$$ 
  
\item $\| (F' \circ f) \cdot f' - (F' \circ g) \cdot g' \|_{L^2}^2$ 
$$\begin{array} {llll}
\leq 2 \| (F' \circ f) \cdot (f' - g') \|_{L^2}^2 + 2 \| ((F' \circ f) - (F' \circ g)) \cdot g' \|_{L^2}^2     \\
\leq 2 \int_0^{+ \infty} (F'(f(x)))^2 (f'- g')^2(x) dx \\
+ 2 \int_0^{+ \infty} (F'(f(x))- F'(g(x)))^2 (g')^2(x) dx   \\
\end{array}$$ 

\noindent Now the first integral is bounded from above by 
$$\left[ \sup_{|y| \leq \| f \|_{L^{\infty}}} |F'(y)|^2 \right] \cdot \| f' - g' \|_{L^2}^2$$ 
\noindent And the second one by 
$$\left[ \sup_{|y| \leq max(\| f \|_{L^{\infty}} ; \| g \|_{L^{\infty}})} |F''(y)|^2 \right] \cdot \|(g)'\|_{L^2}^2 \cdot \| f - g \|_{L^{\infty}}^2 $$      

\noindent Thus, $\| (F' \circ f) \cdot f' - (F' \circ g) \cdot g' \|_{L^2} \leq M'(F,f,g) \cdot \| f - g \|_{H^1((0;+ \infty))}$ with 
$$\begin{array}{ll}
M'(F,f,g) = Const \cdot \max \big( \sup_{|y| \leq \| f \|_{L^{\infty}}} |F'(y)| ; \|g'\|_{L^2} \cdot \\
\hspace{6cm} \left[ \sup_{|y| \leq max(\| f \|_{L^{\infty}} ; \| g \|_{L^{\infty}})} |F''(y)| \right] \big)
\end{array}$$

\end{enumerate}

\noindent In this first step of the proof, we have then proved that, if $F \in C^2(\R)$, \\ 
$\| F \circ f - F \circ g \|_{H^1((0;+ \infty))} \leq D(F,f,g) \| f - g \|_{H^1((0;+ \infty))}$, for any $(f,g)$ in $(C_c^{\infty}(\R) \cap H^2((0;+\infty)))^2$. We have denoted by $D$, the function of $F$, $f$ and $g$ : $D(F,f,g) = Const \cdot (M(F,f,g) + M'(F,f,g))$ where $M(F,f,g)$ and $M'(F,f,g)$ have been defined previously.

\item \it Second step. \rm \\
Using density and the Fatou property of the Hilbert space $H^1((0 ; + \infty))$, we will deduce from the first step that Property $\cc{L} \it ip$ is also satisfied for any $(f,g) \in (H^2((0;+\infty)))^2$ and $F \in C^2(\R,\R)$. \\
\\
\noindent Since $C_c^{\infty}(\R) \cap H^2((0;+\infty))$ is dense into $H^2((0;+\infty))$, there exists a sequel $(f_k)$ included in $C_c^{\infty}(\R) \cap H^2((0;+\infty))$ such that $f_k$ converges to $f$ in $H^2$. Yet $L^{\infty}$ is continuously embedded in $H^1$ so $f_k$ converges to $f$ in $L^{\infty}$ and so, almost everywhere as well. The continuity of $F$ implies that $F \circ f_k$ tends to $F \circ f$ a.e. as $k \longrightarrow +\infty$. \\
\noindent Now, for any $k \in \N$, $F$ and $f_k$ are continuous and, for any bounded interval $\Omega$ of $\R$, the restriction of $f$ to $\Omega$ admits a continuous representative as an $H^2$ function. So $F \circ f_k$ and $F \circ f$ belong to $L^{1,loc}((0;+\infty))$ and, due to Lebesgue's theorem, $F \circ f_k$ converges weakly to $F \circ f$ in $\cc{S'}$ as $k \longrightarrow +\infty$. \\ 
\noindent The same holds for $g$ which can be considered as the limit of a sequel $g_k$. \\
\noindent By assumption, for any $k \in \N$, $F \in C^2(\R)$ and $(f_k, g_k)$ belonging to $(C_c^{\infty}(\R) \cap H^2((0;+\infty)))^2$, then there exists a constant $D(F, f_k,g_k) > 0$ such that $\|F \circ f_k - F \circ g_k \|_{H^1}\leq D(F,f_k,g_k) \cdot \|f_k - g_k \|_{H^1}$ for any $k \in \N$. \\
Since $f_k$ (resp. $g_k$) converges to $f$ (resp. $g$) in $H^2$, $\underline{\lim}_{k \longrightarrow +\infty} \|f_k - g_k \|_{H^1}$ is equal to $\|f - g \|_{H^1}$. Now 
$$D(F, f_k, g_k) = Const \cdot (M(F,f_k,g_k) + M'(F,f_k,g_k))$$
\noindent where 
$$\left \{ \begin{array}{lll}
M(F,f_k,g_k) = \sup_{|y| \leq max(\| f_k \|_{L^{\infty}} ; \| g_k \|_{L^{\infty}})} |F'(y)| \\
M'(F,f_k,g_k) = \sqrt{2} \cdot \max \big( \sup_{|y| \leq \| f_k \|_{L^{\infty}}} |F'(y)| ; \|(g_k)'\|_{L^2}\cdot \\
\hspace{6cm} \left[ \sup_{|y| \leq max(\| f_k \|_{L^{\infty}} ; \| g_k \|_{L^{\infty}})} |F''(y)| \right] \big)
\end{array} \right.$$

\noindent Since $f_k$ converges to $f$ in $L^{\infty}$, $\| f_k \|_{L^{\infty}}$ tends to $\| f \|_{L^{\infty}}$ and so, $M(F,f_k,g_k)$ tends to $M(F,f,g)$ as $k$ tends to $\infty$. On the other hand, $g_k$ converges to $g$ in $H^2$ implies that $(g_k)'$ converges to $g'$ in $L^2$. So $\underline{\lim}_{k \longrightarrow +\infty} M'(F,f_k,g_k) = M'(F,f,g)$. \\
\\
\noindent Thus $\underline{\lim}_{k \longrightarrow +\infty} \|F \circ f_k - F \circ g_k \|_{H^1} \leq D(F,f,g) \cdot \|f - g \|_{H^2}$. \\
Now, the Hilbert space $H^1((0 ; + \infty))$ is such that $\cc{S} \it \subset H^1 \subset \cc{S'}$ and has the Fatou property (cf. \cite{RSi}, p. 15, for the definition of the Fatou property). Thus, $\qqs f \in H^2$, $(F \circ f - F \circ g) \in H^1$ and  
$$\|F \circ f - F \circ g \|_{H^1} \leq Const \cdot D(F,f,g) \cdot \| f - g \|_{H^2}$$ 

\noindent It can also be written : $\|F \circ f - F \circ g \|_{H^1} \leq D'(\| f \|_2, \| g \|_2, \| f' \|_2, \| g' \|_2) \cdot \| f - g \|_{H^2}$ with 
$$\begin{array}{ll}
D'(F,\| f \|_2, \| g \|_2, \| f' \|_2, \| g' \|_2) = Cst \cdot \Big( \sup_{|y| \leq Cst \cdot max(\| f \|_{H^1} ; \| g \|_{H^1})} |F'(y)| \\
+  \max \big( \sup_{|y| \leq Cst \cdot \| f \|_{H^1}} |F'(y)| ; \|g '\|_2 \big) \cdot \left[ \sup_{|y| \leq Cst \cdot max(\| f \|_{H^1} ; \| g \|_{H^1})} |F''(y)| \right] \Big)
\end{array}
$$
\noindent since $H^1$ and $H^2$ are continuously embedded in $L^{\infty}$. \\
To finish with, there lacks to make sure that $D'(F,f,g)$ is a monotone increasing function of its arguments : by definition, $\| f \|_{H^1}$ is a monotone increasing function of $\| f \|_2$ and $\| f' \|_2$. Moreover $M \mapsto \sup_{|y| \leq M} |F'(y)|$ is a monotone increasing function so is $D'$. \fin

\end{itemize}

\end{enumerate}

\begin{rem}~~\\
In the proof, we have also stated that, if $F$ belongs to $C^1(\R)$, then for any $(u,v)$ in $(\prod_{k=1}^2 L^2(\Omega_k))^2$, there exists $C(\| u\|_2, \| v \|_2)$, monotone increasing function of its arguments such that 
$$\| F \circ u - F \circ v \|_2 \leq C(\| u \|_2, \| v \|_2) \cdot \| u - v \|_2$$ 
\noindent which is needed in the second part of Proposition 1 to get the boundedness of the solution and so, its global existence.  
\end{rem}

\begin{cor} ~~\\

\noindent Let $F \in C^2(\R, \R)$ such that $F(0)=0$. If $G : w \mapsto \int_0^w F(\xi) d \xi$ takes only non positive values on $\R$ and if $G \circ u$ belongs to $L^1(\Omega)$ for $u \in L^2(\Omega)$ ($\Omega$ being an interval of $\R$, possibly unbounded), then the unique local solution of Problem $(P)$ guaranteed by the first part of Proposition 1 defines a global solution of $(P)$ i.e. a strongly differentiable $D(\cc A)$-valued function on $[0 ; + \infty)$ that satisfies Problem (P). 
  
\end{cor}

\noindent \bf Proof. \rm \\
Due to Proposition 1, we only have to check that the following sufficient conditions are satisfied :
$$\left \{ \begin{array}{ll}
\| F \circ u - F \circ v \|_2 \leq C(\| u \|_2, \| u_x \|_2, \| v \|_2, \| v_x \|_2) \cdot \| \sqrt{A} u - \sqrt{A} v \|_2 \\
\| \sqrt{A} (F \circ u - F \circ v) \|_2 \leq C(\| u \|_2, \| u_x \|_2, \| v \|_2, \| v_x \|_2) \cdot \| A(u - v) \|_2
\end{array} \right.$$

\noindent for all $(u, v) \in D(A)^2$ and where $C$ is a monotone increasing function of all its arguments. 

\noindent Since the sufficient condition of Proposition 2 is satisfied, it holds for all $(u, v)$ in $D(A)^2$ :
$$\left \{ \begin{array}{ll}
\| F \circ u - F \circ v \|_2 \leq C(\| u \|_2, \| u_x \|_2, \| v \|_2, \| v_x \|_2) \cdot \| u - v \|_{H^1} \\
\| F \circ u - F \circ v \|_{H^1} \leq C(\| u \|_2, \| u_x \|_2, \| v \|_2, \| v_x \|_2) \cdot \| u - v \|_{H^2} 
\end{array} \right.$$

\noindent So the first estimate reads :
$$\begin{array}{lll}
\sum_{k=1}^2 \| F \circ u_k - F \circ v_k \|_2^2 \leq C^2 \left[ \sum_{k=1}^2 \| u_k - v_k \|_2^2 + \sum_{k=1}^2 \| u_{k,x} - v_{k,x} \|_2^2  \right] \\
\leq C^2 \cdot M \cdot \left[ \sum_{k=1}^2 \| u_k - v_k \|_2^2 + \sum_{k=1}^2 a_k \| u_{k,x} - v_{k,x} \|_2^2  \right] \\
\leq C^2 \cdot M \cdot \| \sqrt{A}u - \sqrt{A}v \|_2^2
\end{array}$$
\noindent i.e. $\| F \circ u - F \circ v \|_2 \leq C \cdot \sqrt{M} \cdot \| \sqrt{A}u \|_2$ with $M = \max \left( \dfrac{1}{a_1} ; \dfrac{1}{a_2} ; 1 \right)$. Since $M$ does not depend on $u$ nor on $v$, $C \cdot \sqrt{M}$ is still a monotone increasing function of $\| u \|_2$, $\| u_x \|_2$, $\| v \|_2$ and $\| v_x \|_2$.  \\
Moreover the second estimate given by Proposition 2 reads : \\
$\sum_{k=1}^2 \| F \circ u_k - F \circ v_k \|_2^2 + \| (F \circ u_k)_x - (F \circ v_k)_x \|_2^2$ \\
$\leq C^2 \left[ \sum_{k=1}^2 \| u_k - v_k \|_2^2 + \sum_{k=1}^2 \| u_{k,x} - v_{k,x} \|_2^2 + \sum_{k=1}^2 \| u_{k,xx} - v_{k,xx} \|_2^2  \right]$ \\
Then 
$$\begin{array}{lllllll}
\| \sqrt{A}(F \circ u - F \circ v) \|_2^2   \\
= \| (F \circ u_1 - F \circ v_1) \|_2^2 + a_1 \cdot \| (F \circ u_1 - F \circ v_1)_x \|_2^2 + \| (F \circ u_2 - F \circ v_2) \|_2^2 \\
\hspace{1cm} + a_2 \cdot \| (F \circ u_2 - F \circ v_2)_x \|_2^2 \\
\leq a_2 \cdot \left( \sum_{k=1}^2 \| F \circ u_k - F \circ v_k \|_2^2 + \| (F \circ u_k)_x - (F \circ v_k)_x \|_2^2 \right)\\
\leq C^2 \cdot a_2 \cdot M \cdot \left[ \sum_{k=1}^2 \| u_k - v_k \|_2^2 + \sum_{k=1}^2 \| u_{k,xx} - v_{k,xx} \|_2^2 \right] \\
\hspace{1cm} + C^2 \cdot a_2 \cdot M \cdot \sum_{k=1}^2 \| u_k - v_k \|_2^2 \\
\leq C^2 \cdot a_2 \cdot M \cdot \| Au \|^2 + C^2 \cdot a_2 \cdot M \times 2 \left(\sup_{k \in \{ 1 ; 2 \}} K_k^2 \right) \cdot \| Au - Av \|^2 
\end{array}$$
\noindent where $K_k$ is defined through an Ehrling-Nirenberg-Gagliardo (interpolation) inequality (cf. Adams \cite{ad}, Th 4.14). With $m=p=2$ and $j=1$, it reads for a fixed $\epsilon > 0$ and any $\epsilon_k \leq \epsilon$
$$\| u_{k,x} \|_2^2 \leq K_k \epsilon_k \| u_{k} \|_2^2 + K_k \epsilon_k^{-1} \| u_{k,xx} \|_2^2$$   
\noindent Choosing $\epsilon_k = \epsilon = 1$ leads to the latter estimate. We have then proved that 
$$\| \sqrt{A}(F \circ u - F \circ v) \|_2^2 \leq Const \cdot \| Au - Av \|^2$$
\noindent And the constant is a monotone increasing function of all its arguments : $\| u \|_2$, $\| u_x \|_2$, $\| v \|_2$ and $\| v_x \|_2$.  \fin \\
\\
\noindent \bf Example. \rm If the nonlinearity $F$ is chosen to be defined by $F(u) = - \lambda u^3$ with $\lambda > 0$, it satisfies all the conditions required in Corollary 1 : $F(0)=0$ is clear as well as the regularity ($F \in C^2(\R)$). \\
$G : w \mapsto \int_0^w F(\xi) d \xi$ takes only non positive values on $\R$ : in fact, $G(w) = - \lambda \dfrac{w^4}{4}$ and $w^4 \geq 0$ for any $w \in \R$. At last, since $G$ is continuous and since any $u \in H^2(\Omega)$ is also continuous on $\Omega$, $G \circ u$ lies in $L^{1,loc}(\Omega)$. As for the behaviour at infinity (if $\Omega$ is unbounded), a continuous $L^2$ function is in $L^4$.

\section{Causality of the support of the solution}

\subsection{Spectral solution formula of the linearized problem (LP)}
The aim of this section is to rewrite the spectral solution formula of the linearized problem already studied in \cite{Alreg2} and \cite{Alreg3}, assuming that the initial condition is compactly supported in the first branch $(0 ; + \infty)$. Let us recall that the expression for the unique solution of the problem is given by the following theorem :

\begin{theo} ~~\\
Assume that $f_k : \Omega_k \rightarrow \R$ with $k \in \{ 1 ; 2 \}$, initial data of Problem (LP), described 
in the introduction, are such that $(f_1,f_2) \in D(A)$. 
\\
Then the restriction to $\Omega_1$ of the unique solution of Problem (LP) (such that 
$(u_1,u_2)$ belongs to $C^2(\R^+, \prod_{k=1}^2 
L^2(\Omega_k)) \cap C^1(\R^+, D(\sqrt{A})) 
\cap C^0(\R^+, D(A))$) is given, for $(t,x)$ in $\R^{+*} \times \Omega_1$, by :

$$\begin{array} {llll}
u_1(t,x) =  \frac{1}{2 \pi c^2}  \int_{[\sqrt{a_1};+ \infty)} \cos(\omega t) 
\Im \Big( \frac{2 \omega}{K_1(\omega^2)}   \Big( \int_0^{+\infty} f_1(u)e^{-K_1(\omega^2)(u-x)} 
du \Big) \Big) d\omega \\
\\
 + \frac{1}{2 \pi c^2}  \int_{[\sqrt{a_1};+ \infty)} \cos( \omega t) 
\Im \left( \frac{K_1(\omega^2)-K_2(\omega^2)}{K_1(\omega^2)+K_2(\omega^2)}  \frac{2\omega}{K_1(\omega^2)}  
\left( \int_0^{+\infty} f_1(u)e^{-K_1(\omega^2)(u+x)}du  \right) \right) d\omega\\
\\
 + \frac{1}{2 \pi c^2}  \int_{[\sqrt{a_1};+ \infty)} \cos( \omega t) 
\Im \left( \frac{2K_1(\omega^2)}{K_1(\omega^2)+K_2(\omega^2)}  \frac{2\omega}{K_1(\omega^2)}            
\left( \int_0^{+\infty} f_2(u)e^{-K_2(\omega^2)u-K_1(\omega^2)x}du \right) \right) d\omega 
\end{array}$$

\end{theo}

\noindent $R$, $D(A)$, $D(\sqrt{A})$ are introduced in Section 1.\\
Recall that, for real $\omega$ and $j \in \{ 1 ; 2 \}$

$$K_j(\omega^2) = \left \{ \begin{array} {ll}
\sqrt{\dfrac{a_j - \omega^2}{c^2}} & \mbox{if} \quad \omega^2 \leq 
a_j \\
i \quad  \sqrt{\dfrac{\omega^2 - a_j}{c^2}} & \mbox{if} \quad \omega^2 \geq a_j
\end{array} \right.$$

\noindent Note that the solution is the sum of three terms : the original term, the reflected one which appears due to the discontinuity in the potential and the last one, called "transmitted term", which contains the information on the transmission of the signal from the second branch to the first one. Since our point is to study causality, we will restrict ourselves to a compactly supported initial condition with support in the first branch and the third term will vanish. Rewriting the solution in terms of Fourier transforms, we get :

\begin{theo} ~~\\
Assume that $f_k : \Omega_k \rightarrow \R$ with $k \in \{ 1 ; 2 \}$, initial data of Problem (LP), described 
in the introduction, are such that $(f_1,f_2) \in D(A)$, $f_1$ is compactly supported in $(0; +\infty)$ and $f_2 \equiv 0$. 
\\
Then the restriction to $\Omega_1$ of the unique solution of Problem (LP) (such that 
$(u_1,u_2)$ belongs to $C^2(\R^+, \prod_{k=1}^2 L^2(\Omega_k)) \cap C^1(\R^+, D(\sqrt{A})) 
\cap C^0(\R^+, D(A))$) is given, for $(t,x)$ in $\R^{+*} \times \Omega_1$, by :

$$\begin{array} {ll}
u_1(t,x) = \dfrac{1}{\pi} \cc{F}\it^{-1}_{\omega \mapsto x} \left[\cos(\sqrt{a_1 + c^2 \omega^2}t) \cc{F} \it f_1(\omega) \right]
\\
  \hspace{1.8 cm} + \dfrac{1}{\pi} \cc{F}\it^{-1}_{\omega \mapsto x} \left[\cos(\sqrt{a_1 + c^2 \omega^2}t) \left( \dfrac{\omega c - \sqrt{c^2 \omega^2 -a_2 + a_1}}{\omega c + \sqrt{c^2 \omega^2 -a_2 + a_1}} \right) \cc{F} \it f_1(\omega) \right] 
\end{array}$$

\noindent where the complex square root has been defined such that $\sqrt{c^2 \omega^2 - (a_2 - a_1)}$ =

$$ \left \{ \begin{array} {lll}
\sqrt{c^2 \omega^2 - (a_2 - a_1)} \quad \mbox{if} \quad \omega \geq \sqrt{\frac{a_2 - a_1}{c^2}} \\
i \sqrt{a_2 - a_1 - c^2 \omega^2}  \quad \mbox{if} \quad - \sqrt{\frac{a_2 - a_1}{c^2}} \leq \omega \leq \sqrt{\frac{a_2 - a_1}{c^2}} \\
- \sqrt{c^2 \omega^2 - (a_2 - a_1)} \quad \mbox{if} \quad \omega \leq - \sqrt{\frac{a_2 - a_1}{c^2}}
\end{array} \right.$$ 

\noindent Likewise, for $(t,x)$ in $\R^{+*} \times \Omega_2$, $u_2(t,x) =$  
$$\label{sqrt} \begin{array} {ll}
 - \dfrac{1}{2 \pi} \int_{- \infty}^{\infty} \cos(\sqrt{a_1 + c^2 \omega^2}t) 
e^{i \left( \frac{1}{c} \sqrt{c^2 \omega^2 -a_2 + a_1}  \right) x}
\left( \dfrac{2 \omega c}{\omega c + \sqrt{c^2 \omega^2 -a_2 + a_1}} \right) \cc{F} \it f_1(\omega) d \omega
\end{array}$$

\end{theo}

\noindent The functions of $\omega$ denoted by $C_R$ and $T^{1,2}$ and defined by 
$$C_R(\omega) : = \dfrac{\omega c - \sqrt{c^2 \omega^2 -a_2 + a_1}}{\omega c + \sqrt{c^2 \omega^2 -a_2 + a_1}} \qquad \mbox{and} \qquad T^{1,2}(\omega) := \dfrac{2 \omega c}{\omega c + \sqrt{c^2 \omega^2 -a_2 + a_1}}  $$
are called reflection coefficient and transmission coefficient (from the first branch to the second one) respectively. \\
\\
\bf Proof. \rm \\
\noindent First of all, in the integral defining $u_1$, $\omega$ is changed into $\omega' = \xi_1(\omega^2):=\sqrt{\frac{\omega^2-a_1}{c^2}}$ i.e. $\omega = \sqrt{a_1+c^2 \omega'^2}$ and $d \omega' = \frac{\omega \hspace{0.2cm} d \omega}{c^2 \xi_1(\omega^2)} =  \frac{i}{2c^2} \frac{2 \omega \hspace{0.2cm} d \omega}{K_1(\omega^2)} $. The domain of integration ($(\sqrt{a_1};+ \infty)$) thus becomes $(0;+\infty)$ with the reflection coefficient taking different values depending on the belonging of $\omega$ to $\left[ 0 ; \sqrt{\frac{a_2 -a_1} {c^2}} \right]$ or $\left[\sqrt{\frac{a_2 -a_1} {c^2}} ; + \infty \right)$. \\
\noindent It holds :
$$\begin{array}{lll}
u_1(t,x) =  \dfrac{1}{\pi} \Re \left[ \int_0^{+ \infty} \cos(\sqrt{a_1 + c^2 \omega^2}t) e^{-i \omega x} \cc{F} \it f_1(\omega) d\omega \right]
\\
  \hspace{1 cm} + \dfrac{1}{\pi} \Re \left[\int_0^{\sqrt{\frac{a_2-a_1}{c^2}}} \cos(\sqrt{a_1 + c^2 \omega^2}t) e^{i \omega x} \left( \dfrac{i \omega c - \sqrt{a_2 - a_1 - c^2 \omega^2}}{i \omega c + \sqrt{a_2 - a_1 - c^2 \omega^2}} \right) \cc{F} \it f_1(- \omega)d\omega  \right]\\ 
 \hspace{1 cm} + \dfrac{1}{\pi} \Re \left[\int_{\sqrt{\frac{a_2-a_1}{c^2}}}^{+ \infty} \cos(\sqrt{a_1 + c^2 \omega^2}t) e^{i \omega x}  \left( \dfrac{\omega c - \sqrt{c^2 \omega^2 -a_2 + a_1}}{\omega c + \sqrt{c^2 \omega^2 -a_2 + a_1}} \right) \cc{F} \it f_1(- \omega) d\omega \right]\\ 
\end{array}$$

\noindent In the second step of the proof, $\Re z$ is replaced by $(z + \bar{z})/2$ and $\omega$ is changed into $- \omega$ in the conjugate term. The complex square root of the reflection coefficient is chosen to be analytic in the upper half-plane as in Deutch and Low's paper (cf. \cite{D-L} and the above definition of the complex square root \ref{sqrt}). \\
\noindent As for $u_2$, its expression for a compact supported $f$ is reduced to one term analogous to those in $u_1$ (cf. Theorem 5) : $u_2(t,x) =$ \\
$- \frac{1}{2 \pi c^2}  \int_{[\sqrt{a_1};+ \infty)} \cos( \omega t) 
\Im \left( \frac{2K_2(\omega^2)}{K_1(\omega^2)+K_2(\omega^2)}  \frac{2\omega}{K_2(\omega^2)}            
\left( \int_0^{+\infty} f(u)e^{K_1(\omega^2)u - K_2(\omega^2)x}du \right) \right) d\omega $ \\
\\
\noindent Then the same type of modifications are applied to get the formula of Theorem 5. \fin

\subsection{A Paley/Wiener-type theorem for a function holomorphic in the upper half-plane}

\noindent The aim is to generalize a classical Paley-Wiener Theorem which gives explicit conditions on $L^2(\R)$-functions to be the inverse Fourier transforms of compactly supported functions. This theorem, cited and proved in Rudin (\cite{Ru}), requires analyticity on the whole plane as well as an exponential estimate and a $L^2$-behaviour on the real line. \\
Our solution formula for the linearized problem given in Section 3.1 involves the inverse Fourier transforms of functions which satisfy those requirements except the analyticity on the whole plane. The square roots appearing in both the reflection and transmission coefficients keep them from defining entire functions. They can only be analytic in the upper half-plane (for example) $\Pi_+ = \{ z \in \C / \Im(z)>0 \}$. And we will show that the difference with an entire function is that the support of the inverse Fourier transform of such a function is not compact but only bounded from above.  \\
So our goal here is to adapt the proof of Theorem 19.3 of \cite{Ru} which is based on Cauchy's Theorem. \\
\\
First of all, we need another representation of the square root $\sqrt{c^2 \omega^2 -a_2 + a_1}$ to avoid the branch cuts to lie on the real line (it will be useful to apply Cauchy's Theorem in the proof). Writing, as Ali Mehmeti does in \cite{fam3} 
$$_a \sqrt{z} := \sqrt{|z|} e^{\frac{1}{2}i \arg_a(z)} \quad \mbox{where} \quad z = r e^{i \arg_a(z)}, r \geq 0, \arg_a(z) \in [a ; a + 2 \pi)$$
\noindent we choose $\sqrt{c^2 \omega^2 -a_2 + a_1}$ to be $\left( _{- \pi/2} \sqrt{c \omega - k} \right) \cdot \left( _{- \pi/2} \sqrt{c \omega + k} \right)$ with $k = \sqrt{a_2 - a_1}$. Thus the branch cuts become $\gamma(\pm k) = \{ \pm k - i y /y \geq 0 \}$. This representation does not change the value of the square root on the real line. In particular it is still continuous on the closure of $\Pi_+$ as constructed in the last subsection.

\begin{theo}~~\\
\noindent Let $A$ and $C$ be two positive constants and $g : \C \longrightarrow \C$ be such that :

\begin{itemize}
\item $g$ is analytic everywhere except for the two cuts : \\
$\gamma(\pm k) = \{ \pm k - i y /y \geq 0 \}$ where $k \in \R^{+*}$ 
\item $|g(z)| \leq C e^{A |z|}$, for any $z \in \Pi_+ = \{ z \in \C / \Im(z) > 0 \}$ 
\item $\int_{\R} |g(x)|^2 dx < + \infty$
\end{itemize}

\noindent Define the function G on $\R$ by $G = \cc{F} \it g$ where $\cc{F}$ is the $L^2$-Fourier transform.\\ 
\noindent Then the support of $G$ is a subset of $(-A; + \infty)$ and $G \in L^2(-A;+ \infty)$.
\end{theo}   

\noindent The following proof is an adaptation of Rudin's one (cf. \cite{Ru}). \\
\\
\noindent \bf Proof. \rm \\
\noindent For $\epsilon > 0$ and any real $x$, let $g_{\epsilon}$ be defined by $g_{\epsilon} = g(x) e^{- \epsilon |x|}$. The aim is to show that, for any real $t < -A$ : 
$$\lim_{\epsilon \rightarrow 0} \int_{\R} g_{\epsilon}(x) e^{- itx} dx = 0  \quad (*)$$
Since $\| g - g_{\epsilon} \|_2 \rightarrow 0$ as $\epsilon \rightarrow 0$, Plancherel's Theorem implies that the Fourier transform of $g_{\epsilon}$ (denoted by $G_{\epsilon}$) converges towards the Fourier transform $G$ of the restriction of $g$ to the real axis. Then a subsequence $G_{\epsilon_n}$ converges almost everywhere towards $G$ and, since $lim_{n \rightarrow  +\infty} G_{\epsilon_n}$ vanishes outside $[-A; + \infty)$ (due to $(*)$), so does $G$.    \\
So proving $(*)$ for any real $t < -A$ is sufficient. \\
\\
For any real $\alpha \in [0 ; \pi]$, we define like Rudin :
$$\Gamma_{\alpha}(s) = s e^{i \alpha}, 0 \leq s < + \infty \quad \mbox{and} \quad \Pi_{\alpha} = \{ w / \Re(w e^{i \alpha}) > A \}$$ 
\noindent And for $w \in \Pi_{\alpha}$ :
$$\Phi_{\alpha}(w) = \int_{\Gamma_{\alpha}} g(z) e^{-wz}dz$$
It is clear that $\Phi_{\alpha}$ is defined on $\Pi_{\alpha}$ ($g$ is continuous on the closure of $\Pi_+$ so it is on $\Gamma_{\alpha}$ and $|g(z) e^{- \omega z}| \leq C e^{-s(\Re(\omega e^{i \alpha})-A)}$). It is analytic on $\Pi_{\alpha}$ for $\alpha \in (0;\pi)$ (cf. \cite{Ru}). Moreover it holds for any real $t$ : 
$$\int_{\R} g_{\epsilon}(x) e^{- itx} dx = \Phi_0(\epsilon + it) - \Phi_{\pi} (- \epsilon +it)$$
\noindent The strategy of Rudin is to use Cauchy's Theorem to replace $\Phi_0$ and $\Phi_{\pi}$ by $\Phi_{\pi/2}$ if $t < -A$ and by $\Phi_{-\pi/2}$ for $t>A$. The problem here is that the function $g$ is not entire. We need to avoid the branch cuts to apply Cauchy's Theorem so we define, for $\eta > 0$ : 
$$\Gamma_0^{\eta}= [0 ; k - \eta] \cup [k + \eta ; + \infty) \cup \{k + \eta e^{i \theta}/ \theta \in [0 ; \pi] \} \quad \mbox{and} \quad \Gamma_{\pi}^{\eta}= \{ s / -\bar{s} \in \Gamma_0^{\eta} \}$$
\noindent Then $g$ is analytic on $\Gamma_0^{\eta}$ and on $\Gamma_{\pi}^{\eta}$. \\
Let us also denote by $\Phi_{\alpha}^{\eta}(w) = \int_{\Gamma_{\alpha}^{\eta}} g(z) e^{-wz}dz$ for $\alpha=0$ and $\alpha = \pi$. The interval $[k - \eta ; k + \eta]$ of $\Gamma_0$ has been replaced by a half-circle to bypass the branch point $k$ of the square root so that Rudin's arguments to use Cauchy's Theorem now hold : $\Phi_0^{\eta}(\epsilon + it)$ and $\Phi_{\pi}^{\eta}(\epsilon + it)$ both coincide with $\Phi_{\pi/2}(\epsilon
+ it)$ if $t < -A$. (Note that the domain of $\Phi_{\pi/2}$ is $\Pi_{\pi/2}=\{ z / \Im z < -A \}$.) \\

\begin{lem}~~\\
It holds $\lim_{\eta \rightarrow 0} \Phi_0^{\eta}(w) = \int_{\Gamma_0} g(z) e^{-wz}dz$ and $\lim_{\eta \rightarrow 0} \Phi_{\pi}^{\eta}(w) = \int_{\Gamma_{\pi}} g(z) e^{-wz}dz$. 
\end{lem}   
     
\noindent \bf Proof. \rm \\
The function $g$ is continuous on both half-circles $\gamma_{\eta}(k) = \{k + \eta e^{i \theta}/ \theta \in [0 ; \pi] \}$ and $\gamma_{\eta}(-k) = \{-k + \eta e^{i \theta}/ \theta \in [0 ; \pi] \}$ so 
$$\left|  \int_{\gamma_{\eta}(k)} g(z) e^{-wz}dz \right| \leq \left( \max_{z \in \gamma_{\eta}(k)} |g(z)e^{-wz}| \right) \cdot L(\gamma_{\eta}(k)) $$    
\noindent where $L(\gamma_{\eta}(k))$ is the length of the path $\gamma_{\eta}(k)$ i.e. $\pi \cdot \eta$. Now this path of integration is compact and the integrand is continuous so its maximum is attained for some $z_0$. The half-circle tends to the point $k$ as $\eta$ tends to 0 so $z_0$ tends to $k$ and the maximum value tends to the value of the function at $k$ which is zero. Then $\int_{\gamma_{\eta}(k)} g(z) e^{-wz}dz$ tends to zero as $\eta$ tends to zero. The same result holds for the path $\gamma_{\eta}(-k)$. And the limits of the lemma are a consequence. \fin \\
Thus, for any $t < -A$ : 
$$\begin{array}{lll}
\int_{\R} g_{\epsilon}(x) e^{- itx} dx = \Phi_0(\epsilon + it) - \Phi_{\pi} (- \epsilon +it) \\
\hspace{3cm} = \lim_{\eta \rightarrow 0} \Phi_0^{\eta}(\epsilon + it) - \lim_{\eta \rightarrow 0} \Phi_{\pi}^{\eta}(- \epsilon + it) \\ 
\hspace{3cm} = \Phi_{\pi/2}(\epsilon + it) - \Phi_{\pi/2} (- \epsilon +it)  
\end{array}$$
\noindent It follows $\lim_{\epsilon \rightarrow 0} \int_{\R} g_{\epsilon}(x) e^{- itx} dx = 0$. \fin

\subsection{The main theorem : causality for the nonlinear problem}

\noindent The aim of this last section is to prove that the propagation is causal i.e. that a compactly supported initial signal with vanishing initial velocity is still compactly supported as time goes by. It also means that the propagation of the wave front is not faster than light. Deutch and Low assert in \cite{D-L} that for the linearized problem with a potential barrier, the solution is causally related to the initial condition for propagation to the right. Their idea is that the violation of causality is only apparent since the energy contained in the forward tail of the initial function (a Gaussian centered one) is redistributed. In their situation, the initial condition on $u_t$ does not vanish but is such that their initial wave packet is moving to the right. \\
Our conditions are not so restrictive and our solution is the superposition of right and left-moving waves. Our initial condition is supposed to be compactly supported which is necessary to prove causality properly. The following theorem states that causality still holds when adding a nonlinearity. Our proof is an adaptation of Reed and Simon's Theorem X. 76a in \cite{RS}.     

\begin{theo} ~~\\
Assume that $f$ is such that $(f_1,f_2)$ belongs to $D(A)$ and that the support of $f_1$ is a compact set included in $(0 ; R)$ and $f_2 \equiv 0$. Then the solution of Problem (NLP), given by any of the theorems of the latter section, has the property that the support of $u(t ; \cdot)$ is a subset of
$$\cc{C} \it (\Sigma,t)= \{ x \in \R / \mbox{dist}(x, \Sigma) \leq c t \}$$ 
\end{theo}

\noindent \bf Proof. \rm \\

\begin{lem}~~\\
The reflection and transmission coefficients $C_R(\omega)=\left( \dfrac{\omega c - \sqrt{c^2 \omega^2 -a_2 + a_1}}{\omega c + \sqrt{c^2 \omega^2 -a_2 + a_1}} \right)$ and $T^{1,2}(\omega) = \left( \frac{2 \omega c}{\omega c + \sqrt{c^2 \omega^2 -a_2 + a_1}} \right)$ both define bounded functions of $\omega$ in the upper half-plane $\{ \omega / Im(\omega) > 0 \}$.
\end{lem}

\noindent \bf Proof. \rm \\
Multiplying both coefficients by $\omega c - \sqrt{c^2 \omega^2 -a_2 + a_1}$ leads to 
$$C_R(\omega)= \dfrac{\left( \omega c - \sqrt{c^2 \omega^2 -a_2 + a_1}\right)^2}{a_2 -  a_1} \quad \mbox{and} \quad T^{1,2}(\omega) = \frac{2 \omega c \left( \omega c - \sqrt{c^2 \omega^2 -a_2 + a_1} \right) } {a_2 - a_1}$$ 
\noindent The first expression to study is : $\omega c - \sqrt{c^2 \omega^2 -a_2 + a_1}$. Since it is holomorphic in the upper half-plane, it is bounded provided that its modulus tends to a finite limit when $| \omega |$ tends to $+ \infty$. Denoting by $r$ and $\theta$ the modulus and argument of $\omega$, it holds : $|\omega c - \sqrt{c^2 \omega^2 -a_2 + a_1}|^2 =   $ 
$$\begin{array}{lll}
\left( cr \cos(\theta) - M(r) \cos(A(\theta)) \right)^2 + \left( cr \sin(\theta) - M(r) \sin(A(\theta)) \right)^2 \\
\mbox{with} \quad M(r,\theta):= \left( \left( c^2 r^2 \cos^2(2 \theta) - a_2 + a_1 \right)^2 + c^4 r^4 \sin^2(2 \theta) \right)^{1/4}\\
\mbox{and} \quad A(r, \theta):= \frac{1}{2} \arctan \left( \frac{c^2 r^2 \sin(2 \theta)}{c^2 r^2 \cos(2 \theta)- a_2 + a_1 } \right) 
\end{array}$$
\noindent Now an easy computation gives $M(r,\theta):= cr \left[ 1 + \dfrac{(a_1-a_2) \cos(2 \theta)}{2 c^2 r^2} + o \left( \dfrac{1}{r^2} \right) \right]$ and $A(r, \theta) \equiv \theta$ so that $|\omega c - \sqrt{c^2 \omega^2 -a_2 + a_1}| = \dfrac{(a_2-a_1) \cos(2 \theta)}{2 cr} + o \left( \dfrac{1}{r} \right)$. \\
Thus $\lim_{r \rightarrow \infty} C_R(r e^{i \theta}) = 0$ and $\lim_{r \rightarrow \infty} T^{1,2}(r e^{i \theta}) = \cos(2 \theta)$. \fin \\
\noindent Note that if $\omega$ is real, the reflection coefficient tends to 0 and the transmission one to 1 as the frequency $\omega$ tends to infinity which is a reasonable and unsurprising behaviour. \\
\\
\noindent Let us first consider $u_1$ : \\
\noindent The proof of Reed and Simon (cf. \cite{RS}, p. 309) is adaptable. Only the first part of the proof concerning the linear part has to be changed in the following way. The solution reads in terms of Fourier transforms :
$$\begin{array} {ll}
u_1(t,x) =  \cc{F}\it^{-1}_{\omega \mapsto x} \left[\cos(\sqrt{a_1 + c^2 \omega^2}t) \cc{F} \it f_1(\omega) \right]
\\
  \hspace{1.8 cm} + \cc{F}\it^{-1}_{\omega \mapsto x} \left[\cos(\sqrt{a_1 + c^2 \omega^2}t) \left( \dfrac{\omega c - \sqrt{c^2 \omega^2 -a_2 + a_1}}{\omega c + \sqrt{c^2 \omega^2 -a_2 + a_1}} \right) \cc{F} \it f_1(\omega) \right] 
\end{array}$$
\noindent Let us check the assumptions of Theorem 6 for $g$ defined by 
$$g(z)= \left[\cos(\sqrt{a_1 + c^2 z^2}t) \left( \dfrac{cz - \sqrt{c^2 z^2 -a_2 + a_1}}{cz + \sqrt{c^2 z^2 -a_2 + a_1}} \right) \cc{F} \it f_1(z) \right]$$ 

\noindent First of all, since $f_1$ is an $L^2$-function with compact support contained in $(0 ; R)$, the Fourier transform of $f_1$ is an entire analytic function and there exists a real constant $C_1$ such that 
$$|\cc{F} \it f_1(\omega)| \leq C_1 e^{|\omega| R}$$
\noindent (cf. \cite{Ru}, introduction of Chapter 19 for this classical result). \\
\noindent Further since the roots drop out in the power series, $\cos(\sqrt{a_1 + c^2 \omega^2}t)$ is also entire and satisfy, for some real constant $C_2$ 
$$\left| \cos(\sqrt{a_1 + c^2 \omega^2}t) \right| \leq C_2 e^{|\omega| ct}$$
\noindent Now, since the modulus of the reflection coefficient is bounded (cf. Lemma 2), there is a real constant $C$ such that 
$$|g(\omega)| \leq C e^{|\omega| (R+ct)}$$
\noindent i.e. the estimate of Theorem 6 is satisfied with $A:=R + ct$. \\
\noindent Furthermore, since $f_1$ is an $L^2$-function, then its Fourier transform lies in $L^2(\R)$ as well. And since the reflection coefficient is bounded on $\R$ as well as the cosinus, $g$ belongs to $\L^2(\R)$.\\  
Since the reflection coefficient is an analytic function in the upper half-plane only and it has two branch cuts : $\gamma(\pm k)$ with $k = \sqrt{a_2 - a_1}$, the function $G$ defined on $\R$ by $G(x) = \int_{\R} g(\omega) e^{-i \omega x} d\omega$ is such that the support of $G$ is a subset of $(-A;+\infty)$ by Theorem 6. Then that of $x \mapsto G(-x)$ is included in $(-\infty;A)$ i.e. $(-\infty;R+ct)$. \\
Since $u_1(t, \cdot)$ is the restriction to $\R^+$ of this function, it is compactly supported in $(0;R + ct)$ i.e. its support is in $\cc{C} \it ((-R ; R),ct)$. The rest of the proof remains unchanged. \\
\\
As for $u_2$, it is rewritten as : $u_2(t,x) =$ \\
$ - \frac{1}{2 \pi c^2}  \int_{- \infty}^{+ \infty} \cos(\sqrt{a_1 + c^2 \omega^2}t) 
e^{i \omega x} e^{i \left( - \omega - \frac{1}{c} \sqrt{c^2 \omega^2 -a_2 + a_1}  \right) x}
\left( \frac{2 \omega c}{\omega c + \sqrt{c^2 \omega^2 -a_2 + a_1}} \right) \cc{F} \it f_1(\omega) d \omega$

\noindent Now, since the transmission coefficient and the additional factor $e^{i \left( - \omega - \frac{1}{c} \sqrt{c^2 \omega^2 -a_2 + a_1}  \right) x}$ are bounded (cf. Lemma 2), the same arguments as before give the compactness of the support of $u_2$ as regards with the space variable $x$. \fin


\newpage

\begin{thebibliography}{99} 

\bibitem{ad} R. Adams. \it Sobolev spaces. \rm Acad. Press, 1975. \\
\bibitem{fam-3} F. Ali Mehmeti. \it Problèmes de transmission pour des
équations des ondes linéaires et quasilinéaires. \rm Séminaire Equations aux
Dérivées Partielles Hyperboliques et Holomorphes (1982/83), J. Vaillant (ed.),
Travaux en cours, Hermann, Paris, p. 75-96, 1984. \\
\bibitem{fam-2} F. Ali Mehmeti. \it A characterization of a
generalized $C^{\infty}$-notion on nets. \rm Integral Equations and
Operator Theory, \bf 9\rm, p. 753-766, 1986. \\
\bibitem{fam-1} F. Ali Mehmeti. \it Regular Solutions of Transmission
and Interaction Problems for Wave Equations. \rm Math. Meth. in the
Appl. Sci., \bf 11\rm, p. 665-685, 1989.\\
\bibitem{fam1} F. Ali Mehmeti. \it
Nonlinear Waves in Networks. \rm Mathematical Research, vol. 80, Akademie Verlag, Berlin, 1994.\\
\bibitem{fam2} F. Ali Mehmeti. \it Spectral Theory and $L^{\infty}$ - time
Decay Estimates for Klein-Gordon Equations on Two Half Axes with Transmission
: the Tunnel Effect. \rm Math. Meth. in the Appl. Sci., \bf 17\rm, p. 697-752,
1994.\\ 
\bibitem{fam3} F. Ali Mehmeti. \it Transient Waves in Semi-Infinite
Structures : the Tunnel Effect and the Sommerfeld Problem. \rm Mathematical Research, vol. 91, Akademie Verlag, Berlin, 1996.\\
\bibitem{Alreg1} F. Ali Mehmeti, V. Régnier. \it Splitting of energy
of dispersive waves in a star-shaped network. \rm Z. Angew. Math. Mech., \bf83\rm, No 2, p. 105-118, 2003.\\
\bibitem{Alreg2} F. Ali Mehmeti, V. Régnier. \it Réflexion retardée pour des paquets d'ondes dispersives sur un réseau en forme d'étoile. \rm C. R. Acad. Sci. Paris, Ser. I \bf337 \rm p. 645-648, 2003.\\
\bibitem{Alreg3} F. Ali Mehmeti, V. Régnier. \it Delayed reflection of the energy
flow at a potential step for dispersive wave packets. \rm Math. Meth. Appl. Sci., \bf27\rm, p. 1145-1195, 2004.\\
\bibitem{Lum} F. Ali Mehmeti, J. von Below and S. Nicaise, edit. \it
Partial differential equations on multistructures, Lecture Notes in Pure and
Applied Mathematics\rm, Marcel Dekker, 2001. \\
\bibitem{App} J. Appell. \it Nonlinear superposition operators. \rm Cambridge Tracts in Mathematics, 95, Cambridge University Press, 1990. \\
\bibitem{Br} H. Brézis. \it Analyse fonctionnelle, théorie et applications. \rm 
Collection Ma\-thé\-ma\-ti\-ques Appliquées pour la Maîtrise, Masson, Paris, 1983. \\
\bibitem{CD} E. Croc, Y. Dermenjian. \it
Spectral analysis of a multistratified acoustic strip Part II : asymptotic behaviour of solutions for a simple stratification. \rm SIAM Journal on Mathematical Analysis \bf
27\rm(6), p. 1632-1652, 1996. \\
\bibitem{yas} Y. Daikh. \it Temps de passage de paquets d'ondes de basses fréquences ou limités en bandes de fréquences par une barrière de potentiel. \rm Thèse de doctorat, Valenciennes, France, 2004. \\
\bibitem{D-L} J. M. Deutch, F. E. Low. \it
Barrier Penetration and Superluminal Velocity. \rm Annals of Physics \bf
228\rm, p. 184-202, 1993. \\
\bibitem{Em} T. Emig. \it Propagation of an
electromagnetic pulse through a waveguide with a barrier : an exact solution
within classical electrodynamics. \rm Diplomarbeit, Institut f\"ur
Theoretische Physik, K\"oln, Germany,1996. \\ 
\bibitem{Nim} A. Enders, G. Nimtz. \it On superluminal barrier traversal. \rm
J. Phys. I France, \bf2\rm, p. 1693-1698, 1992.
\bibitem{Ha} A. Haraux. \it Nonlinear Evolution Equations, Global Behavior of Solutions. \rm Lecture Notes in Mathematics, Springer-Verlag, Berlin, Heidelberg, New York, 317 p., 1980. \\
\bibitem{Hla} J. Hladik, M. Chrysos. \it Introduction à la mécanique
quantique. \rm Dunod, 212 p., 2000. \\
\bibitem{Jr} K. J\"orgens. \it Das Anfangswertproblem im Grossen f\"ur eine Klasse nichtlinearer Wellengleichungen. \rm Math. Z., \bf77\rm, p.295-308, 1961. \\
\bibitem{Kl} S. Klainermann. \it Global existence for Nonlinear Wave Equations. \rm Comm. Pure Appl. Math., \bf 33\rm, p.
43-101, 1980.\\ 
\bibitem{K} S. B. Kuksin. \it On Squeezing and flow of energy for nonlinear
wave equations. \rm Geometric and Functional Analysis, \bf 5\rm, No 4, p.
668-701, 1995.\\ 
\bibitem{MSW} B. Marshall, W. Strauss, S. Wainger. \it $L^p$ - $L^q$
Estimates for the Klein-Gordon Equation. \rm J. Math. Pures et Appl., \bf 59\rm, p. 417-440, 1980. \\
\bibitem{Mi} K. Mihalin\u{c}i\'{c}. \it Time Decay Estimates for the Wave
Equation with Transmission and Boundary Conditions. \rm Dissertation.
Technische Universit\"{a}t Darmstadt, Germany, 1998. \\
\bibitem{nith} S. Nicaise. \it Diffusion sur les espaces ramifiés. \rm
Thèse de doctorat, Mons, Belgique, 1986.\\
\bibitem{ni2} S. Nicaise. \it Spectre des réseaux topologiques finis. \rm
Bull. Sci. Math. (2) \bf 111\rm, No 4, p. 401-413, 1987. \\
\bibitem{RS} M. Reed, B. Simon. \it Methods of modern mathematical physics, II : Fourier analysis, self-adjointness. \rm Academic Press, New York, San Francisco, London, 361 p., 1975. \\
\bibitem{reg} V. Régnier. \it Delayed reflection in a stratified acoustic strip. \rm Math. Meth. Appl. Sci., \bf28\rm, p. 185-203, 2005.\\
\bibitem{Rott} K. Rottbrand. \it Time-Dependent Plane Wave Diffraction by a
Half-Plane : Explicit Solution for Rawlins' Mixed Initial Boundary Value
Problem\rm. Z. Angew. Math. Mech., \bf78\rm, No 5, p. 321-334, 1998. \\  
\bibitem{Ru} W. Rudin. \it Real and Complex Analysis. \rm TMH Edition, 1966.\\ 
\bibitem{RSi} T. Runst and W. Sickel. \it Sobolev Spaces of Fractional Order, Nemytskij Operators and Nonlinear Partial Differential Equations. \rm Walter de Gruyter and Co, Berlin, 1952.\\ 
\bibitem{S} I. Segal. \it Nonlinear semi-groups. \rm Annals of Math. \bf 78\rm, No 2, p. 339-364, 1962. \\
\bibitem{Sh} J. Shatah. \it Global Existence of Small Solutions to Nonlinear Evolution Equations. \rm J. Diff. Eq. \bf 46\rm, No 3, p. 409-425, 1982. \\
\bibitem{St} W. Strauss. \it Nonlinear wave equations. \rm CBMS \bf 73\rm, AMS, Providence, Rhode Island, 1989. \\
\bibitem{Weid} J. Weidmann. \it Spectral theory of ordinary differential operators. \rm Lecture Notes in Mathematics, 1258, Springer-Verlag, New York, 1987. \\



\end{thebibliography}
\end{document}